\documentclass [final,3p,times]{elsarticle-preprint}

\usepackage{amssymb,amsmath,color}

\newtheorem{theorem}{Theorem}

\newtheorem{proposition}[theorem]{Proposition}
\newtheorem{corollary}[theorem]{Corollary}

\newcommand{\R}{{\mathbb R}}
\newcommand{\N}{{\mathbb N}}
\newcommand{\be}[1]{\begin{equation}\label{#1}}
\newcommand{\ee}{\end{equation}}
\newcommand{\eqn}[1]{\eqref{#1}}
\newcommand{\ird}[1]{\int_{\R^d}{#1}\,dx}
\newcommand{\irdmu}[1]{\int_{\R^d}{#1}\,d\mu}
\newcommand{\nrm}[2]{\|{#1}\|_{L^{#2}(\R^d)}}
\newcommand{\nrmu}[2]{\|{#1}\|_{L^{#2}(\R^d,\,d\mu)}}
\renewcommand{\(}{\left(}
\renewcommand{\)}{\right)}
\newcommand{\OU}{\textsf N}

\newenvironment{proof}{\noindent\emph{Proof.\/} }{\unskip\null\hfill$\square$\vskip 0.3cm}

\journal{}

\begin{document}

\begin{frontmatter}

\title{Improved intermediate asymptotics for the heat equation\tnoteref{Thanks,Copyright}}
\tnotetext[Thanks]{This work has been partially supported by the \emph{Fondation Sciences Math\'ematiques de Paris} and by the ANR projects \emph{CBDif-Fr} and \emph{Evol}.}
\tnotetext[Copyright]{\copyright\,2009 by the authors. This paper may be reproduced, in its entirety, for non-commercial purposes.}

\author[Ceremade]{Jean-Philippe Bartier}
\ead{bartier@ceremade.dauphine.fr}
\author[Toulouse]{Adrien Blanchet}
\ead{adrien.blanchet@univ-tlse1.fr}
\author[Ceremade]{Jean Dolbeault}
\ead{dolbeaul@ceremade.dauphine.fr}
\author[Bilbao]{Miguel Escobedo}
\ead{miguel.escobedo@ehu.es}
\address[Ceremade]{CEREMADE (UMR CNRS no. 7534), Universit\'e Paris-Dauphine, Place de Lattre de Tassigny, 75775 Paris C\'edex~16, France.}
\address[Toulouse]{GREMAQ (UMR CNRS no. 5604 and INRA no. 1291), Universit\'e de Toulouse, 
21 all\'ee de Brienne, 31000 Toulouse, France.}
\address[Bilbao]{Departamento de Matem\'aticas, 
Universidad del Pa\'{\i}s Vasco, Barrio Sarriena s/n, 48940 Lejona (Vizcaya), Spain.}
\date{June 19, 2009}
\begin{abstract}
This letter is devoted to results on intermediate asymptotics for the heat equation. We study the convergence towards a stationary solution in self-similar variables. By assuming the equality of some moments of the initial data and of the stationary solution, we get improved convergence rates using entropy / entropy-production methods. We establish the equivalence of the exponential decay of the entropies with new, improved functional inequalities in restricted classes of functions. This letter is the counterpart in a linear framework of a recent work on fast diffusion equations, see~\cite{BDGV}. Results extend to the case of a Fokker-Planck equation with a general confining potential.\end{abstract}


\end{frontmatter}

Consider the \emph{heat equation} in the euclidean space,
\begin{equation}\label{Eqn:H}
\frac{\partial u}{\partial t}=\Delta u\quad t>0\,,\;x\in\R^d
\end{equation}
with an initial condition $u_0\in L^1(\R^d)$. By writing $u=u_+-u_-$ where $u_+$ and $u_-$ are respectively the positive and negative parts of $u$ and solving \eqref{Eqn:H} with initial data $(u_0)_+$ and $(u_0)_-$, we may reduce the problem to the case of a nonnegative function, corresponding to a nonnegative initial condition $u_0$, without restriction. The heat equation being linear, we can assume without loss of generality that $u_0$ is a probability measure so that in the sequel of this note $\int_{\R^d}u_0\,dx=1=\int_{\R^d}u(t,x)\,dx$ for any $t\ge 0$. Getting decay rates and even an asymptotic expansion for large values of $t$ is completely standard, see for instance \cite{MR1183805}. However, a few details and some notations will be useful for later purpose.

\medskip First of all, as a straightforward consequence of the expression of the Green function, $G(t,x,y):=(4\pi t)^{-d/2}\,e^{-\frac{|x-y|^2}{4t}}$, any solution $u$ of \eqref{Eqn:H} can be written as $u(t,x)=\int_{\R^d}u_0(y)\,G(t,x,y)\,dy$ and therefore uniformly decays like $O(t^{-d/2})$ since, as $t\to\infty$, $u(t,x)\sim G(t,x,0)$. It is also classical to estimate the decay of $u(t,\cdot)- G(t,\cdot,0)$ in various $L^p(\R^d)$ norms. Such estimates are called \emph{intermediate asymptotics} estimates. The point is to determine the first term of an asymptotic expansion of the solution as $t\to\infty$. For instance, as we shall see below, it can be proved that $\nrm{u(t,\cdot)-G(t,\cdot,0)}1=O(t^{-1/2})$ as $t\to\infty$.

\medskip The \emph{entropy method} can be used among various other approaches to obtain such an estimate. It relies on the logarithmic Sobolev inequality and goes as follows. First consider the time-dependent rescaling
\be{Eqn:rescaling}
u(t,x)=R^{-d}\,v\(\log R,x/{R}\)\quad\mbox{with}\quad R=R(t):=\sqrt{1+2\,t}\,,\quad t>0\,,\;x\in\R^d\,.
\ee
If $u$ is a solution of \eqn{Eqn:H}, then $v$ solves the \emph{Fokker-Planck equation}
\begin{equation}\label{Eqn:FP}
\frac{\partial v}{\partial t}=\Delta v+\nabla\cdot(x\,v)
\end{equation}
with same initial condition $v(t=0,\cdot)=u_0$. Let $v_\infty(x):=(2\pi)^{-d/2}\,e^{-|x|^2/2}$ be the unique stationary solution of~\eqn{Eqn:FP} with mass $1$, and define $d\mu:=v_\infty\,dx$ as the Gaussian measure. We denote by $L^p(\R^d)$ and $L^p(\R^d,\,d\mu)$ the Lebesgue spaces corresponding respectively to Lebesgue's measure and to the Gaussian measure. Understanding the intermediate asymptotics for $u$ amounts to study the convergence of $v$ to~$v_\infty$, as $t\to\infty$. Define the \emph{entropy} by $\mathcal E_1[w]:=\int_{\R^d} {w\log w}\,d\mu$. Let $v$ be a solution of~\eqn{Eqn:FP} and define $w(t,\cdot):=v(t,\cdot)/v_\infty$, $w_0:=w(t=0,\cdot)$. Then $\frac d{dt}\,\mathcal E_1[w(t,\cdot)]=-\mathcal I\!_1[w(t,\cdot)]$ where $\mathcal I\!_1$ is the \emph{Fisher information} defined by $\mathcal I\!_1[w]:=\irdmu{w\,|\nabla \log w|^2}$. Gross' \emph{logarithmic Sobolev inequality} exactly amounts to $\mathcal E_1[v/v_\infty]\leq\frac 12\,\mathcal I\!_1[v/v_\infty]$ and so, it follows that
\begin{equation*}
 \mathcal E_1[w(t,\cdot)]\le \mathcal E_1[w_0]\,e^{-2\,t}\quad\forall\,t\ge 0\,.
\end{equation*}
By the \emph{Csisz\'ar-Kullback inequality}, see for instance \cite{MR1801751}, we get $\nrm{v(t,\cdot)-v_\infty}1^2\le\frac 1{4}\,\mathcal E_1[w(t,\cdot)]$ and deduce that
\[
\nrm{v(t,\cdot)-v_\infty}1\le\frac 12\,\sqrt{\mathcal E_1[w_0]}\,e^{-t}\quad\forall\,t\ge 0\,.
\]
Undoing the change of variables~\eqref{Eqn:rescaling} and observing that $u_\infty(t,x):=R(t)^{-d}\,v_\infty\(x/R(t)\)=G(t+1/2,\cdot,0)$, we finally get
\[
\nrm{u(t,\cdot)-u_\infty(t,\cdot)}1\le \frac 12\,\sqrt{\frac{\mathcal E_1[w_0]}{1+2\,t}}\quad\forall\;t\ge 0\;,
\]
which establishes the claimed estimate, namely: $\nrm{u(t,\cdot)-G(t,x,0)}1\le O\(t^{-1/2}\)$ as $t\to\infty$. Such an estimate is quite classical. The above method is known as the \emph{Bakry-Emery method} or \emph{entropy / entropy-production method} and also provides a proof of the logarithmic Sobolev inequality. See \cite{MR1447044,Arnold-Markowich-Toscani-Unterreiter01} for some references on this topic, in the context of partial differential equations. 

\medskip By combining $L^1(\R^d)$ and $L^\infty(\R^d)$ estimates using H\"older's inequality, we get that 
\[
\nrm{u(t,\cdot)-G(t,\cdot,0)}p\le O\Big(t^{-\frac 1{2\,p}\,(1+(p-1)\,d)}\Big)\quad\mbox{as}\quad t\to\infty\;.
\]
In a $L^2(\R^d)$ framework, a much more detailed description can be achieved using a spectral decomposition. If $v$ is a solution of \eqn{Eqn:FP}, then $w=v/v_\infty$ is a solution of the \emph{Ornstein-Uhlenbeck equation}
\be{Eqn:OU}
\frac{\partial w}{\partial t}=\Delta w-x\cdot\nabla w
\ee
with initial data $w_0=u_0/v_\infty$. Notice that $\irdmu{w_0}=1$ and, as a consequence, $\irdmu{w(t,\cdot)}=1$ for all $t\geq 0$. Define by $(H_k)_{k\in\N^d}$ the sequence of Hermite type polynomials (see for instance \cite{weisstein:hp}) acting on $x=(x_1,x_2\ldots x_d)\in\R^d$, such that $H_k(x):=\prod_{j=1}^d h_{k_j}(x_j)$ where $h_n(y):=(-1)^n\,(n!)^{-1/2}\,e^{y^2/2}\frac{d^n}{dy^n}\big(e^{-y^2/2}\big)$, $y\in\R$ and $k=(k_1,...,k_d) \in \N^d$. These functions provide an orthonormal family of eigenfunctions in $L^2(\R^d,\,d\mu)$ which spans the eigenspaces of the Ornstein-Uhlenbeck operator, that is $-\(\Delta H_k-x\cdot\nabla H_k\)=|k|\,H_k$, where $|k|:=\sum_{j=1}^dk_j$. Up to a scaling, $(h_n)_{n\in\N}$ is the usual family of Hermite polynomials on $\R$. 

If $w_0$ satisfies the orthogonality condition
\be{Cond:Orthogonality}
\irdmu{w_0\,H_k}=0\quad \forall\,k\in\N^d\;\mbox{such that}\;0<|k|<n\,,
\ee
then an improved rate of convergence follows, in the sense that
\[
\left\|w(t,\cdot)-1\right\|_{L^2(\R^d,\,d\mu)}\le e^{-n\,t}\left\|w_0-1\right\|_{L^2(\R^d,\,d\mu)}\quad\forall\,t\ge 0\,.
\]
If \eqref{Cond:Orthogonality} initially holds, we indeed have $\irdmu{w(t,\cdot)\,H_k}=0$ for any $t\ge 0$ and any $k\in\N^d$such that $0<|k|<n$. Then, since $\frac d{dt}\nrmu{w(t,\cdot)-1}2^2=-2\irdmu{\left|\nabla w(t,\cdot)\right|^2}$, the conclusion holds using the following result.
\begin{proposition}[Improved Poincar\'e inequality]\label{Prop:improvedPoincare} Assume that $w\in L^2(\R^d)$ is such that $\irdmu w=1$ and the condition $\irdmu{w\,H_k}=0$ holds for any $k\in\N^d$ such that $0<|k|<n$. Then the following inequality holds, with optimal constant:
\[
\nrmu{w-1}2^2\leq\frac 1n\,\nrmu{\nabla w}2^2\,.
\]\end{proposition}
The proof is no more than a straightforward rewriting of the Rayleigh quotient $\nrmu{\nabla w}2^2/\nrmu{w-1}2^2$ under the appropriate orthogonality condition. Notice that polynomials $H_k$ are of degree $|k|$ so that the Condition \eqref{Cond:Orthogonality} can be rephrased in terms of moment conditions. See~\cite{MR1183805,kim2009hoa} for further results in this direction.

\medskip It is natural to search for improved estimates of convergence also in $L^p(\R^d)$ with $p\in[1,2)$ by looking for improved functional inequalities whenever condition \eqref{Cond:Orthogonality} is fulfilled. We may for instance quote \cite{MR2409469} in which improvements on the constant, but not on the rates, have been achieved for $p=1$. 

For any $p\in(1,2]$, consider the \emph{generalized entropy}
\[
\mathcal E_p[w]:=\irdmu{\frac {w^p-1}{p-1}}\,.
\]
This definition is consistent with the definition of $\mathcal E_1$ because, under the condition $\irdmu {\kern -1.5pt w}=1$, $\mathcal E_p[w]=\irdmu{\kern -1.5pt\frac {w^p-w}{p-1}}\to\mathcal E_1[w]$ as $p\to 1$. The functional $\mathcal E_p$ controls the convergence in $L^p(\R^d,\,d\mu)$ using a generalized Csisz\'ar-Kullback inequality. In \cite{MR1951784,BDIK}, it has been proved that $\|w-1\|_{L^p(\R^d,\,d\mu)}^2\le \frac 1p\,{2^{2/p}}\,\max\big\{\,\|w\|_{L^p(\R^d,\,d\mu)}^{2-p},1\big\}\,\mathcal E_p[w]$, for any $p\in [1,2]$. Since $\nrmu w1=1$, we have $1\le\nrmu wp^p=1+(p-1)\,\mathcal E_p[w]$, and so
\be{eq:ck}
\|w-1\|_{L^p(\R^d,\,d\mu)}\le \mathcal A_p\left(\mathcal E_p[w]\right)\quad\mbox{with}\quad\mathcal A_p(s):=\frac{2^{1/p}}{\sqrt{p}}\,\Big[1+(p-1)\,s\Big]^{1-p/2}\sqrt{s}\;.
\ee
Next, assume that $\irdmu{w\,H_k}=0$ for any $k\in\N^d$ such that $0<|k|<n$ and consider the \emph{generalized Poincar\'e inequalities,} with $p\in[1,2]$, namely
\be{Beckner}
\mathcal E_p[w]\le\mathcal B_{n,p}\,\irdmu{\left|\nabla w^{p/2}\,\right|^2}\quad\forall\;w\in H^1(\R^d,\,d\mu)\,.
\ee
Such inequalities have been established for $n=1$ by W. Beckner in \cite{MR954373} with optimal constant $\mathcal B_{1,p}=2/p$ for the Gaussian measure. By the same method, it has been shown in \cite{0528} that for a larger class of measures $d\mu$, if \eqref{Beckner} holds for $p=1$ and $p=2$, for some positive constants $\mathcal B_{n,1}$ and $\mathcal B_{n,2}$ respectively, then it also holds for any $p\in(1,2)$ with
\begin{equation}
 \label{bnp}
\textstyle \mathcal B_{n,p}=\frac 1{p-1}\, \left[1-\((2-p)/p\)^{\,\mathcal B_{n,1}/(2\,\mathcal B_{n,2})}\,\right]\,\mathcal B_{n,2}\;.
\end{equation}
By the logarithmic Sobolev inequality and the improved Poincar\'e inequality, see Proposition~\ref{Prop:improvedPoincare}, we know that $\mathcal B_{n,1}\le 2$ and $\mathcal B_{n,2}=1/n$. Hence it follows that $\mathcal B_{n,p}\le\frac 1{p-1}\, \left[1-\((2-p)/p\)^n\,\right]\,\frac 1n$. On the other hand, as in \cite{Arnold-Markowich-Toscani-Unterreiter01}, if $w$ is a solution of \eqref{Eqn:OU}, then
\be{E-EP}
\frac d{dt}\,\mathcal E_p[w(t,\cdot)]=-\frac 4p\irdmu{\big|\,\nabla w^{\,p/2}\,\big|^2}\,.
\ee
If \eqref{Cond:Orthogonality} is satisfied, we conclude using \eqref{Beckner} and \eqref{eq:ck} that any solution of~\eqref{Eqn:OU} with initial data $w_0$ satisfies 
\[
\mathcal E_p[w(t,\cdot)]\le \mathcal E_p[w_0]\,e^{-2\,\lambda(n,p)\,t}\quad\mbox{and}\quad\|w(t,\cdot)-1\|_{L^p(\R^d,\,d\mu)}\le \mathcal A_p\left(\mathcal E_p[w_0]\right)\,e^{-\lambda(n,p)\,t}\quad\forall\;t\ge0\,,
\]
with $\lambda(n,p):=\frac 2p\,n\,(p-1)\left[1-\((2-p)/p\)^n\right]^{-1}$. The last estimate holds because, for any $t \ge 0$,
\begin{equation*}
 \|w(t,\cdot)-1\|_{L^p(\R^d,\,d\mu)} \le \mathcal A_p\left(\mathcal E_p[w(t,\cdot)]\right) \le \mathcal A_p\left(\mathcal E_p[w_0]\,e^{-2\,\lambda(n,p)\,t} \right) \le \mathcal A_p\left(\mathcal E_p[w_0]\right)e^{-\lambda(n,p)\,t}
\end{equation*}
Notice that $\lambda(1,p)=1$ and $\lambda(n,2)=n$. Nothing is gained as $p\to 1$, since $\lim_{p\to 1}\lambda(n,p)=1$ is independent of $n$.

On the other hand, by H\"older's inequality, we have for free that $\nrmu{w-1}p\le\nrmu{w-1}2$. Hence, if $w$ is a solution of~\eqref{Eqn:OU} with initial data $w_0$, we know that $\nrmu{w(t,\cdot)-1}p\le e^{-n\,t}\,\nrmu{w_0-1}2$ as $t\to\infty$, for any $p\in[1,2]$, if \eqref{Cond:Orthogonality} is satisfied. By interpolation, we recover the rates of~\cite{MR1183805,kim2009hoa}. However, this is not satisfactory since neither $\nrmu{w_0\!-\!1}p$ nor $\mathcal E_p[w_0]$ are involved in the right hand side of the above estimate.

\medskip Consider first the case $p=1$. An alternative approach is suggested by the method of~\cite{BBDGV-Cras,BBDGV}, which applies to the fast diffusion equation $\frac{\partial u}{\partial t}=\Delta u^m$ for $m<1$. By assuming some uniform bound on the initial data, which is preserved along the evolution, it is possible to relate the asymptotic rate for intermediate asymptotics with the spectrum of the linearized operator. We can indeed observe that $\nrmu{w_0-1}2^2\le\nrmu{w_0-1}1\,\nrmu{w_0-1}\infty\le\frac 12\sqrt{\mathcal E_1[w_0]}\,\nrmu{w_0-1}\infty$ using H\"older's inequality and the Csisz\'ar-Kullback inequality. This proves that
\[
\nrmu{w(t,\cdot)-1}1^2\le \tfrac 12\,\nrmu{w_0-1}\infty\,\sqrt{\mathcal E_1[w_0]}\,e^{-n\,t}\quad\mbox{as}\;t\to\infty
\]
if \eqref{Cond:Orthogonality} is satisfied initially. Still, this provides neither an estimate of $\irdmu{w(t,\cdot)\,\log w(t,\cdot)}$ nor a functional inequality which improves upon the logarithmic Sobolev inequality. To prove such an inequality, we keep following the strategy of \cite{BBDGV}. A simple but key idea is to observe that the functions defined for any $p\in[1,2]$ by $h_p(0)=1$, $h_p(1)=p/2$ and, for any $s \in (0,1)\cup(1,\infty)$ by $h_p(s):=[{s^p-1-p\,(s-1)}]/[(p-1)\,|s-1|^2]$ if $p>1$, $h_1(s):=[{s\,\log s-(s-1)}]/{|s-1|^2}$, are continuous, nonnegative, decreasing on $\R^+$ and achieve their maximum at $0$. Define on $L^\infty(\R^d)$ the functional 
\[
\mathcal H_p[w]:=\nrm w\infty^{2-p}\sup_{x\in\R^d}h_p(w(x))=\nrm w\infty^{2-p}\;h_p\(\,\inf_{x\in\R^d}w(x)\)\,.
\]
\begin{theorem}[Improved logarithmic Sobolev inequality]\label{THM:LSIimprovedlogsob} Assume that $w\in L^\infty_+(\R^d)$ is such that $\irdmu w=1$ and satisfies the condition $\irdmu{w\,H_k}=0$ for any $k\in\N^d$ such that $0<|k|<n$. Then the following inequality holds, with optimal constant:
\[
\irdmu{w\,\log w}\leq\frac{\mathcal H_1[w]}n\irdmu{\frac{|\nabla w|^2}w}\,.
\]\end{theorem}
\begin{proof} We may indeed observe that by the Poincar\'e inequality and using the definition of $\mathcal H_1$, we get
\[
\irdmu{\frac{\vert\nabla w\vert^2}{w}}\ge\frac 1{\nrm w\infty}\irdmu{\vert\nabla w\vert^2}\ge \frac n{\nrm w\infty}\irdmu{|w-1|^2}\ge\frac n{\mathcal H_1[w]}\,\irdmu{w\log w }\,.
\]
The optimality of the constant can be checked by a lengthy but elementary computation using the functions $w_\varepsilon^k:=H_k(x)\,\chi\big(x\,\varepsilon^{1/(2n)}\big)+C_\varepsilon^k$ for some smooth truncation function $\chi$ such that $0\le\chi\le 1$, $\chi\!\equiv \!1$ on $B(0,1)$ and $\chi\equiv 0$ in $\R^d\setminus B(0,2)$. Here for $k\in\N^d$ is such that $|k|=n$ and the constant $C_\varepsilon^k$ is chosen so that $\irdmu{w_\varepsilon^k}=1$.\end{proof}


As a consequence of the Maximum Principle applied to the heat equation~\eqref{Eqn:H} and the fact that to $u_0=v_\infty$ corresponds a self-similar solution of~\eqref{Eqn:H}, namely $u(t,x)=G(t+\tfrac 12,x,0)$, we have the estimate
\[
\mathcal H_1[w(t,\cdot)]\le\mathcal H_1[w_0]\quad\forall\;t\ge 0\,.
\]
By applying Theorem~\ref{THM:LSIimprovedlogsob}, we obtain a new result of decay for $\mathcal E_1[w(t,\cdot)]$ with a constant which is exactly $\mathcal E_1[w_0]$, to the price of a rate which is less than $2\,n$.
\begin{corollary}[Improved decay rate of the entropy]\label{Cor:improvedRate} Let $w$ be a solution of~\eqref{Eqn:OU} with a nonnegative bounded initial data $w_0\in L^1(\R^d,\,d\mu)$ such that $\irdmu {w_0}=1$ and \eqref{Cond:Orthogonality} is satisfied. Then
\[
\mathcal E_1[w(t,\cdot)]\le\mathcal E_1[w_0]\,e^{-{n\,t}/\,{\mathcal H_1[w_0]}}\quad\forall\;t\ge 0\,.
\]
\end{corollary}
This result is actually equivalent to Theorem~\ref{THM:LSIimprovedlogsob}, as follows by differentiating the above inequality at $t=0$ (for which equality is trivially satisfied) and using the fact that $-\irdmu{|\nabla w_0|^2/w_0}=\frac d{dt}\,\mathcal E_1[w(t,\cdot)]_{|t=0}\le\mathcal E_1[w_0]\,\frac d{dt}\,{e^{-{n\,t}/\,{\mathcal H_1[w_0]}}}_{|t=0} $. What we have achieved is a global, improved exponential decay of the entropy $\mathcal E_1$ in a restricted class of functions. To simplify even further, for any $\varepsilon\in(0,1)$ and $n\in\N^*$, consider the set $\mathcal X_\varepsilon^n:=\{w\in L^1(\R^d,\,d\mu)\,:\,1-\varepsilon\le w\le1+\varepsilon\;\mbox{\sl a.e. and the condition $\irdmu{w\,H_k}=0$ holds for any $k\in\N^d$ such that $0<|k|<n$}\,\}$, which is appropriate to handle the optimality case corresponding to $\varepsilon\to 0_+$. The best constant in Theorem~\ref{THM:LSIimprovedlogsob} is indeed asymptotically equivalent to the sharp rate of convergence in Corollary~\ref{Cor:improvedRate}, in the sense that $\lim_{\varepsilon\to 0_+}\inf_{w\in\mathcal X_\varepsilon^n}n/{\mathcal H_1[w]}=\lim_{\varepsilon\to 0_+}{n}/[{(1+\varepsilon)\,h(1-\varepsilon)}]=2\,n$.

\medskip For simplicity, we have considered only the case $p=1$, but the method also applies to any $p\in(1,2)$. We obtain an improved version of~\eqref{Beckner} under the restriction that $w\in L^1(\R^d,\,d\mu)$ is bounded nonnegative and the condition $\irdmu{w\,H_k}=0$ holds for any $k\in\N^d$ such that $0<|k|<n$. With $\mathcal B_{n,1}=4\,\mathcal H_1[w]/n$ and $\mathcal B_{n,2}=1/n$, we get $\mathcal B_{n,p}\le\mathcal K[n,p,w]:=(n\,(p-1))^{-1}\,\left[1-\((2-p)/p\)^{2\,\mathcal H_1[w]}\right]$ by~\eqref{bnp}. Using the entropy / entropy-production identity~\eqref{E-EP}, the fact that $\mathcal K[n,p,w(t,\cdot)]\le \mathcal K[n,p,w_0]$ and the generalized Csisz\'ar-Kullback inequality~\eqref{eq:ck}, we obtain
\be{Eqn:pCase1}
\mathcal E_p[w(t,\cdot)]\le\mathcal E_p[w_0]\,e^{-\frac{4\,t}{p\,\mathcal K[n,p,w_0]}}\quad\mbox{and}\quad\nrmu{w-1}p\le\mathcal A_p\left(\mathcal E_p[w_0]\right)\,e^{-\frac{2\,t}{p\,\mathcal K[n,p,w_0]}}\quad\forall\;t\ge 0\,.
\ee

Alternatively, an elementary computation as in the proof of Theorem~\ref{THM:LSIimprovedlogsob} gives a similar result:
\[
\frac 4{p^2}\irdmu{\big|\,\nabla w^{p/2}\,\big|^2}=\irdmu{w^{p-2}\,\vert\nabla w\vert^2}\ge\frac 1{\nrm w\infty^{2-p}}\irdmu{\vert\nabla w\vert^2}\ge \frac n{\nrm w\infty^{2-p}}\irdmu{|w-1|^2}\ge\frac n{\mathcal H_p[w]}\,\mathcal E_p[w]
\]
if  $\irdmu w=1$ and the condition $\irdmu{w\,H_k}=0$ holds for any $k\in\N^d$ such that $0<|k|<n$. This proves that
\[
\mathcal E_p[w]\le\frac 4{p^2}\,\frac{\mathcal H_p[w]}n\irdmu{\big|\,\nabla w^{p/2}\,\big|^2}\,.
\]
Using~\eqref{E-EP} and~\eqref{eq:ck}, this proves that any solution of \eqref{Eqn:OU} with initial data in $w_0\in L^1\cap L^\infty(\R^d,\,d\mu)$ satisfies
\be{GenEntDecay-GenNrmDecay}
\mathcal E_p[w(t,\cdot)]\le\mathcal E_p[w_0]\,e^{-{n\,p\,t}/\,{\mathcal H_p[w_0]}}\quad\mbox{and}\quad\|w-1\|_{L^p(\R^d,\,d\mu)}\le\mathcal A_p\left(\mathcal E_p[w_0]\right)\,\,e^{-{n\,p\,t}/{(2\,\mathcal H_p[w_0])}}\quad\forall\;t\ge 0\,.
\ee
Comparing the rates of \eqref{Eqn:pCase1} and \eqref{GenEntDecay-GenNrmDecay} is a natural question. In the limit $\varepsilon\to 0$, $\inf_{w\in \mathcal X_\varepsilon^n}{\mathcal H_p[w]}\sim\sup_{w\in \mathcal X_\varepsilon^n}{\mathcal H_p[w]}\to p/2$ and it follows that $\lim_{\varepsilon\to 0}\frac 4{p\,\mathcal K[n,p,w_0]}=\frac 4p\,n\,(p-1)/[1-((2-p)/p)^p]<2\,n=\lim_{\varepsilon\to 0}\frac{n\,p\,t}{\mathcal H_p[w_0]}$. Hence, at least in the regime $\varepsilon\to 0$, \eqref{GenEntDecay-GenNrmDecay} is a better estimate in terms of rates than \eqref{Eqn:pCase1}. Undoing the change of variables~\eqref{Eqn:rescaling}, we have achieved a detailed result on improved $u_0$.
\begin{corollary}[Improved intermediate asymptotics for the heat equation] Let $p\in[1,2]$ and assume that $u_0$ is a probability measure such that $w_0=u_0/v_\infty$ is bounded and satisfies the condition $\ird{u_0\,H_k}=0$ for any $k\in\N^d$ such that $0<|k|<n$. If $u$ is the solution of~\eqref{Eqn:H} with initial condition $u_0$, then
\[
\nrm{u(t,\cdot)-u_\infty(t,\cdot)}p\le (2\pi)^{-\frac d2\,(1-\frac 1p)}\,\mathcal A_p\left(\mathcal E_p[w_0]\right)\,\,(1+2\,t)^{-\frac{n\,p}{4\,\mathcal H_p[w_0]}-\frac d2\,(1-\frac 1p)}\quad\forall\;t\ge 0\,.
\]
\end{corollary}
The proof relies on the remark that $\nrm{u(t,\cdot)-u_\infty(t,\cdot)}p\le \nrm{u_\infty(t,\cdot)}\infty^{1-\frac 1p}\,\nrmu{w(t,\cdot)-1}p$ where $u_\infty(t,\cdot):=G(t+1/2,\cdot,0)$. The conclusion holds using $\nrm{u_\infty(t,\cdot)}\infty=(2\pi\,R^2)^{-d/2}$ with $R=\sqrt{1+2\,t}$. 

\medskip Up to now, we have considered the simple case of the harmonic potential, $V(x)=\frac 12\,|x|^2$. As in~\cite{0528}, the previous results can be extended to more general potentials as follows. Consider $V\in W^{1,2}_{\rm loc}\cap W^{2,2}_{\rm loc}(\R^d)$ such that \hbox{$\int_{\R^d}e^{-V(x)}dx=1$}, and define the probability measure $d\mu(x):={e^{-V(x)}}dx$ in~$\R^d$, which generalizes the Gaussian measure. Under the above conditions on $V$, the logarithmic Sobolev inequality holds (resp. \eqref{Beckner} for $p=1$) for some positive constant (resp. for $\mathcal B_{1,1}>0$). The Ornstein-Uhlenbeck operator $\OU:=-\Delta+\nabla V\cdot\nabla$ is essentially self-adjoint on $L^2(d\mu)$, has a non-degenerate eigenvalue $\lambda_0=0$ and a spectral gap $\lambda_1>0$. According to~\cite[Theorem~2.1]{MR1736202}, $\OU$ has a pure point spectrum without accumulation points. Since $\lim_{k\to\infty}\lambda_k=\infty$, then by~\cite[Theorem~XIII.64]{MR0493421}, the eigenfunctions of $\OU$ form a complete basis of $L^2(\R^d,\,d\mu)$. We shall denote the eigenvalues by $\lambda_k$, $k\in\N$, and by $E_k$ the corresponding eigenspaces. 

Theorem~\ref{THM:LSIimprovedlogsob} adapts without changes. Assume that $w\in L^\infty_+(\R^d)$ is such that $\irdmu w=1$. Then
\[
\irdmu{w\,\log w}\leq\frac{\mathcal H_1[w]}{\lambda_n}\irdmu{\frac{|\nabla w|^2}w}
\]
under the orthogonality condition: $w \in \left( \bigcup_{k=1}^{n-1}E_k\right)^\bot$, that is $\irdmu{w\,f_k}=0$ for any $f_k\in E_k$, $k=1$, $2$,\ldots$n-1$. Next, consider the solution $w$ of the Ornstein-Uhlenbeck equation
\be{eqn:OUgeneral}
\frac{\partial w}{\partial t}=-\OU\,w=\Delta w - \nabla V\cdot\nabla w\,,
\ee
with initial condition $w_0\in\left( \bigcup_{k=1}^{n-1}E_k\right)^\bot\cap L^\infty(\R^d)$ is such that $\irdmu{w_0}=1$. With the same definition as above for $\mathcal E_p$, for any solution of \eqref{eqn:OUgeneral} with initial data $w_0$, \eqref{GenEntDecay-GenNrmDecay} is now replaced by
\[
\mathcal E_p[w(t,\cdot)]\le\mathcal E_p[w_0]\,e^{-{\lambda_n\,p\,t}/\,{\mathcal H_p[w_0]}}\quad\mbox{and}\quad\|w-1\|_{L^p(\R^d,\,d\mu)}\le\mathcal A_p\left(\mathcal E_p[w_0]\right)\,\,e^{-{\lambda_n\,p\,t}/\,{(2\,\mathcal H_p[w_0])}}\quad\forall\;t\ge 0\,.
\]

\medskip Let us conclude this letter by some comments and open questions. It is standard in entropy / entropy-production methods that determining sharp rates of convergence in an evolution equation is equivalent to finding sharp constants in functional inequalities, as we have seen in the case of the heat equation: the rate of convergence in $L^2(\R^d,\,d\mu)$ is given by the Poincar\'e inequality, while the rate of convergence in entropy, which controls the $L^1(\R^d,\,d\mu)$ norm, is related with the logarithmic Sobolev inequality. This is also true for nonlinear diffusion equations, see for instance~\cite{MR1940370}. In this case, a breakthrough came from the observation that uniform norms can also be used, see \cite{MR2224869,BBDGV-Cras,BBDGV}, to the price of a restricted functional framework. This allows to relate nonlinear quantities of entropy type with spectral properties of the linearized problem, in an appropriate functional space and, again, to relate sharp rates with best constants, see~\cite{BDGV}. As long as nonlinear evolution problems are concerned, only a few invariant quantities are usually available: the mass and the position of the center of mass of the solution, for instance. In linear evolution problems, we can impose an arbitrary number of orthogonality conditions, which are preserved along the evolution. Improved rates of convergence are then expected, even when measured with nonlinear quantities like the entropy. Various attempts have been done, see for instance \cite{MR2409469}, but the question has been left open for many years. Such ideas have been partially explored by R.J. McCann, including in the linear case (see \cite{CE-MC}), based on considerations on an appropriate Hessian matrix. Our approach provides a simpler and elementary answer under restrictions which are natural in view of \cite{BBDGV}. It also raises a number of questions concerning the optimality of the new functional inequalities from a variational point of view, the convergence of minimizing sequences and the symmetry of the eventual minimizers.


\end{document}